# Milnor numbers of projective hypersurfaces with isolated singularities

June Huh

ABSTRACT

Let $V$ be a projective hypersurface of fixed degree and dimension which has only isolated singular points. We show that, if the sum of the Milnor numbers at the singular points of $V$ is large, then $V$ cannot have a point of large multiplicity, unless $V$ is a cone. As an application, we give an affirmative answer to a conjecture of Dimca and Papadima.

## 1. Main results

**1.1**
Let $h$ be a nonzero homogeneous polynomial of degree $d \geqslant 1$ in the polynomial ring $\mathbb{C}[z_0, \ldots, z_n]$. In order to avoid trivialities, we assume throughout that $n \geqslant 2$. We write $V(h)$ for the projective hypersurface $\{h = 0\} \subseteq \mathbb{P}^n$. Associated to $h$ is the *gradient map* obtained from the partial derivatives

$$\mathrm{grad}(h) : \mathbb{P}^n \dashrightarrow \mathbb{P}^n, \qquad z \longmapsto \left(\frac{\partial h}{\partial z_0} : \cdots : \frac{\partial h}{\partial z_n}\right).$$

The *polar degree* of $h$ is the degree of the gradient map of $h$. The polar degree of $h$ depends only on the set $\{h = 0\}$ [DP03]. If $V(h)$ has only isolated singular points, then the polar degree is given by the formula

$$\deg\bigl(\mathrm{grad}(h)\bigr) = (d-1)^n - \sum_{p \in V(h)} \mu^{(n)}(p),$$

where $\mu^{(n)}(p)$ is the Milnor number of $V(h)$ at $p$. See [DP03, Section 3], and also [FM12, Huh12a].

THEOREM 1. *Suppose $V(h)$ has only isolated singular points, and let $m$ be the multiplicity of $V(h)$ at one of its points $x$. Then*

$$\deg\bigl(grad(h)\bigr) \geqslant (m-1)^{n-1},$$

*unless $V(h)$ is a cone with the apex $x$.*

Equivalently, we have

$$(d-1)^n \geqslant (m-1)^{n-1} + \sum_{p \in V(h)} \mu^{(n)}(p),$$

unless $V(h)$ is a cone with the apex at $x$. Therefore, for projective hypersurfaces of given degree and dimension, the existence of a point of large multiplicity should be compensated by a smaller sum of the Milnor numbers at the singular points, unless the hypersurface is a cone.

 The author was partially supported by NSF grant DMS-0943832.

It is interesting to observe how badly the inequality fails when $V(h)$ is a cone over a smooth hypersurface in $\mathbb{P}^{n-1} \subseteq \mathbb{P}^n$. In this case, the polar degree is zero, but the apex of the cone has multiplicity $d$. Both the multiplicity and the sum of the Milnor numbers are simultaneously as large as possible with respect to the degree and the dimension.

The inequality also crucially depends on the assumption that $V(h)$ has only isolated singular points. This can be most clearly seen from Gordan-Noether counterexamples to Hesse's claim that the polar degree is zero if and only if the hypersurface is a cone [Hes51, Hes59, GN76]. For example, consider the threefold in $\mathbb{P}^4$ defined by

$$h = z_3^{d-1} z_0 + z_3^{d-2} z_4 z_1 + z_4^{d-1} z_2, \qquad d \geqslant 3.$$

In this case, the degree of the gradient map of $h$ is zero, $V(h)$ has a point of multiplicity of $d-1$, but $V(h)$ is not a cone.

**1.2**

Let $(V, \mathbf{0})$ be the germ of an isolated hypersurface singularity at the origin of $\mathbb{C}^n$. Associated to the germ are the *sectional Milnor numbers* introduced by Teissier [Tei73]. The $i$-th sectional Milnor number of the germ, denoted $\mu^{(i)}$, is the Milnor number of the intersection of $V$ with a general $i$-dimensional plane passing through $\mathbf{0}$.

We obtain Theorem 1 from the following refinement.

THEOREM 2. *Suppose $V(h)$ has only isolated singular points, and let $\mu^{(n-1)}$ be the $(n-1)$-th sectional Milnor number of $V(h)$ at one of its points $x$. Then*

$$\deg\bigl(\mathrm{grad}(h)\bigr) \geqslant \mu^{(n-1)},$$

*unless $V(h)$ is a cone with the apex $x$.*

The Minkowski inequality for mixed multiplicities says that the sectional Milnor numbers always form a log-convex sequence [Tei77]. In other words, we have

$$\frac{\mu^{(n)}}{\mu^{(n-1)}} \geqslant \frac{\mu^{(n-1)}}{\mu^{(n-2)}} \geqslant \cdots \geqslant \frac{\mu^{(i)}}{\mu^{(i-1)}} \geqslant \cdots \geqslant \frac{\mu^{(1)}}{\mu^{(0)}}.$$

Since $\mu^{(0)}$ is one and $\mu^{(1)}$ is one less than the multiplicity at the point, we see that Theorem 2 indeed implies Theorem 1.

The inequality of Theorem 2 is tight relative to the degree and the dimension. In Proposition 17 we show that, for each $d \geqslant 3$ and $n \geqslant 2$, there is a degree $d$ hypersurface in $\mathbb{P}^n$ with one singular point, for which the equality holds in Theorem 2.

**1.3**

Our interest in Theorems 1 and 2 arose from the study of homalodial polynomials. A *homaloidal polynomial* is a homogeneous polynomial whose gradient map is a birational transformation of $\mathbb{P}^n$. A *homaloidal hypersurface* is the projective hypersurface defined by a homaloidal polynomial. See [CRS08] for a motivated introduction.

Dolgachev showed in [Dol00] that there are exactly three homaloidal plane curves, up to a linear change of homogeneous coordinates:

– a nonsingular conic $h = z_0^2 + z_1^2 + z_2^2 = 0$.
– the union of three nonconcurrent lines $h = z_0 z_1 z_2 = 0$.

- the union of a conic and one of its tangent $h = z_0(z_1^2 + z_0 z_2) = 0$.

In contrast, there are abundant examples of homaloidal hypersurfaces in $\mathbb{P}^n$ when $n \geqslant 3$.

- Any relative invariant of a regular prehomogeneous space is homalodal [ES89, EKP02].
- Projective duals of certain scroll surfaces are homaloidal [CRS08].
- Determinants of generic sub-Hankel matrices are homaloidal [CRS08].
- The union of a cone over a homaloidal hypersurface in $\mathbb{P}^{n-1} \subseteq \mathbb{P}^n$ and a general hyperplane is homaloidal [FM12].
- There are infinitely many polytopes such that almost all polynomials having any one of them as the Newton polytope are homaloidal [Huh12b].

In particular, the last construction shows that there are irreducible homaloidal hypersurfaces of any given degree $d \geqslant 3$ in the projective space of dimension $n \geqslant 3$. Dimca and Papadima conjectured in [DP03, Section 3] that none of them has only isolated singular points.

CONJECTURE 3. *There are no homaloidal hypersurfaces of degree $d \geqslant 3$ with only isolated singular points in the projective space of dimension $n \geqslant 3$.*

We use Theorem 2 to give an affirmative answer to this conjecture.

THEOREM 4. *A projective hypersurface with only isolated singular points has polar degree 1 if and only if it is one of the following, after a linear change of homogeneous coordinates:*

- $(n \geqslant 2, d = 2)$ *a smooth quadric*
$$h = z_0^2 + \cdots + z_n^2 = 0.$$
- $(n = 2, d = 3)$ *the union of three nonconcurrent lines*
$$h = z_0 z_1 z_2 = 0.$$
- $(n = 2, d = 3)$ *the union of a smooth conic and one of its tangent*
$$h = z_0(z_1^2 + z_0 z_2) = 0.$$

Theorem 4 shows that, for projective hypersurfaces of given degree and dimension, the sequence of possible values for the sum of the Milnor numbers necessarily contains a gap, except for quadric hypersurfaces and cubic plane curves. Similar, but stronger results concerning the sum of Tjurina numbers can be found in the works of du Plessis and Wall [dPW99, dPW01]. Other important evidence in support of Conjecture 3 was provided by [Dim01, CRS08, Ahm10].

We close this introduction by posing the problem of finding other forbidden values for the sum of the Milnor numbers at the singular points of a degree $d$ hypersurface in $\mathbb{P}^n$, for general $d$ and $n$. See Conjecture 20.

**1.4**

We now provide a brief overview of the paper.

In Section 2, we formulate and prove a Lefschetz hyperplane theorem with an assigned base point for projective hypersurface complements. The main argument involves

- a pencil of hyperplane sections which has only isolated singular points with respect to a Whitney stratification [Tib02a, Tib02b], and

- a generalized Zariski theorem on the fundamental groups of plane curve complements [Dim92, Section 4.3].

In Section 3, we prove Theorem 2, using the Lefschetz hyperplane theorem of the previous section. We provide an example showing that the inequality of Theorem 2 is sharp.

In Section 4, we use Theorem 2 to show that all the singularities of a homaloidal hypersurface with only isolated singular points are necessarily simple of type $A$. The proof of Conjecture 3 is then obtained from the results of du Plessis and Wall, and Dimca [dPW01, Dim01]. We close with a brief discussion of projective hypersurfaces with polar degree 2.

Acknowledgements

The author thanks Nero Budur, Igor Dolgachev, Anatoly Libgober, and Mircea Mustaţă for helpful discussions. He thanks Radu Laza for providing the proof of Conjecture 20 for cubic threefolds.

## 2. Lefschetz theorem with an assigned base point

Let $D(h)$ be the hypersurface complement $\{h \neq 0\} \subseteq \mathbb{P}^n$. Hamm's Lefschetz theory shows that, if $H$ is a general hyperplane in $\mathbb{P}^n$, then

$$\pi_i\big(D(h), D(h) \cap H\big) = 0 \quad \text{for} \quad i < n.$$

See [Ham83, HL85]. The purpose of this section is to refine this result by allowing hyperplanes to have an assigned base point.

THEOREM 5. *If $H_x$ is a general hyperplane passing through a point $x$ in $\mathbb{P}^n$, then*

$$\pi_i\big(D(h), D(h) \cap H_x\big) = 0 \quad \text{for} \quad i < n,$$

*unless*

(i) *one of the components of $V(h)$ is a cone with the apex $x$, or*

(ii) *the singular locus of $V(h)$ contains a line passing through $x$.*

Since $D(h)$ and $D(h) \cap H_x$ are homotopic to CW-complexes of dimensions $n$ and $n-1$ respectively, the vanishing of the homotopy groups implies

$$H_i\big(D(h), D(h) \cap H_x\big) = 0 \quad \text{for} \quad i \neq n.$$

*Example* 6. Let $V(h)$ be the plane curve consisting of a nonsingular conic containing $x$, the tangent line to the conic at $x$, and a general line passing through $x$. Then

$$H_1\big(D(h), D(h) \cap H_x\big) \simeq H_1\big(S^1 \times S^1, S^1\big) \simeq \mathbb{Z}.$$

*Example* 7. Let $V(h)$ be the cone over a smooth hypersurface of degree $d$ in $\mathbb{P}^{n-1} \subseteq \mathbb{P}^n$ with the apex $x$. Then

$$H_{n-1}\big(D(h), D(h) \cap H_x\big) \simeq H_{n-1}\big(D(h) \cap \mathbb{P}^{n-1}, D(h) \cap H_x \cap \mathbb{P}^{n-1}\big) \simeq \mathbb{Z}^{(d-1)^{n-1}}.$$

It seems reasonable to expect that the condition on the singular locus of $V(h)$ is also necessary. However, the author does not know this.

The rest of this section is devoted to the proof of Theorem 5.

**2.1**

We start with a characterization of the apex of irreducible cones in $\mathbb{P}^n$.

LEMMA 8. *Let $V$ be a subvariety of positive dimension $k+1$ in $\mathbb{P}^n$. Then the following conditions are equivalent for a point $x$ in $\mathbb{P}^n$.*

  (i) *$V$ is a cone with the apex $x$.*
 (ii) *For any point $y$ of $V$ different from $x$, the line joining $x$ and $y$ is contained in $V$.*
(iii) *If $E_x$ is a general codimension $k$ linear subspace in $\mathbb{P}^n$ containing $x$, then every irreducible component of $V \cap E_x$ is a line containing $x$.*
(iv) *If $E_x$ is a general codimension $k$ linear subspace in $\mathbb{P}^n$ containing $x$, then some irreducible component of $V \cap E_x$ is a line containing $x$.*

The irreducibility assumption is clearly necessary in order to deduce (iii) from (iv).

*Proof.* The equivalence of the first three conditions is standard, and (iii) implies (iv). We show that (iv) implies (ii), using a pointed version of the Fano variety of lines in $V$.

Consider the Grassmannians

$$G_1 := \{L_x \mid L_x \text{ is a line in } \mathbb{P}^n \text{ containing } x\} \simeq \mathrm{Gr}(1,n),$$
$$G_2 := \{E_x \mid E_x \text{ is a codimension } k \text{ linear subspace in } \mathbb{P}^n \text{ containing } x\} \simeq \mathrm{Gr}(n-k,n),$$

and the incidence correspondence

$$\mathscr{I} := \{(L_x, E_x) \mid L_x \subseteq V \cap E_x\} \subseteq G_1 \times G_2.$$

Let $\mathrm{pr}_1, \mathrm{pr}_2$ be the projections from $\mathscr{I}$

$$\begin{array}{ccc} & \mathscr{I} & \\ \mathrm{pr}_1 \swarrow & & \searrow \mathrm{pr}_2 \\ G_1 & & G_2. \end{array}$$

We compute the dimension of the image of $\mathrm{pr}_1$, the variety of lines through $x$ contained in $V$. Our assumption (iv) says that $\mathrm{pr}_2$ is generically surjective. Since $\mathrm{pr}_2$ is generically finite in general, it follows that the dimension of $\mathscr{I}$ is equal to that of $G_2$. Note also that

$$\mathrm{pr}_1^{-1}(L_x) \simeq \begin{cases} \mathrm{Gr}(n-k-1, n-1) & \text{if } L_x \subseteq V, \\ \varnothing & \text{if } L_x \not\subseteq V. \end{cases}$$

Therefore

$$\dim \mathrm{Im}(\mathrm{pr}_1) = \dim \mathrm{Gr}(n-k,n) - \dim \mathrm{Gr}(n-k-1, n-1) = k.$$

Next, consider the incidence correspondence

$$\mathcal{I} := \{(L_x, p) \mid p \in L_x\} \subseteq \mathrm{Im}(\mathrm{pr}_1) \times V,$$

and the associated projections

$$\begin{array}{ccc} & \mathcal{I} & \\ \pi_1 \swarrow & & \searrow \pi_2 \\ \mathrm{Im}(\mathrm{pr}_1) & & V. \end{array}$$

$\pi_1$ is a bundle of projective lines, and therefore the dimension of $\mathcal{I}$ is $k+1$. $\pi_2$ is injective over the open subset $\{p \neq x\}$, because there is at most one line containing $p$ and $x$ which is contained in $V$. Since $V$ is assumed to be irreducible, the previous two sentences imply that $\pi_2$ is surjective. In other words, for any point $y$ of $V$ different from $x$, the line joining $x$ and $y$ is contained in $V$. □

**2.2**

Let $X$ be a smooth projective variety of dimension $n$, and let $A$ be a general codimension 2 linear subspace of a fixed ambient projective space of $X$. One of the conclusions of the classical Lefschetz theory is the isomorphism

$$H_{i+1}(X, X_c) \simeq H_{i-1}(X_c, X_c \cap A), \qquad i < n-1,$$

where $X_c$ is a general member of the pencil of hyperplane sections of $X$ associated to $A$ [Lam81, Section 3.6]. By induction, one has the vanishing

$$H_i(X, X_c) = 0, \qquad i < n.$$

The aim of this subsection is to state a generalization due to Tibăr [Tib02a, Tib02b]. We state these results in the generality that we need, not necessarily in the generality of the original papers. See also [Tib07, Section 10.1].

We work in the following setting:

- $V$ is a closed subset of a projective variety $Y$.
- $X$ is the quasi-projective variety $Y \setminus V$.
- $\mathscr{W}$ is a Whitney stratification of $Y$ such that $V$ is a union of strata.
- $A$ is a codimension 2 linear subspace of a fixed ambient projective space of $Y$.
- $\mathscr{W}|_{Y \setminus A}$ is the Whitney stratification of $Y \setminus A$ obtained by restricting $\mathscr{W}$.
- $\mathscr{P}_A$ is the pencil of hyperplanes containing the axis $A$. We write

$$\pi : Y \setminus A \longrightarrow \mathscr{P}_A$$

  for the map sending $p$ to the member of $\mathscr{P}_A$ containing $p$.
- $\mathbb{Y}$ is the blow-up of $Y$ along $Y \cap A$. We write

$$p : \mathbb{Y} \longrightarrow \mathscr{P}_A$$

  for the map which agrees with $\pi$ on $Y \setminus A$.
- $\mathscr{S}$ is a Whitney stratification of $\mathbb{Y}$ which extends $\mathscr{W}|_{Y \setminus A}$.

By a Whitney stratification we mean a complex analytic partition which satisfies the Whitney regularity conditions and the frontier condition. For generalities on Whitney stratifications we refer to [GWPL76, LT10] and references therein.

DEFINITION 9. The *singular locus* of $p$ with respect to $\mathscr{S}$ is the following closed subset of $\mathbb{Y}$:

$$\operatorname{Sing}_{\mathscr{S}} p := \bigcup_{\mathcal{S} \in \mathscr{S}} \operatorname{Sing} p|_{\mathcal{S}}.$$

We say that $\mathscr{P}_A$ has only *only isolated singular points* with respect to $\mathscr{S}$ if $\dim \operatorname{Sing}_{\mathscr{S}} p \leqslant 0$.

The singular locus of $p$ is a closed subset of $\mathbb{Y}$ because $\mathscr{S}$ is a Whitney stratification. The notion of isolated singular points in this generalized sense has proved its value, for example, in the works of Lê [Lê87, Lê92].

We are now ready to introduce the theorem of Tibăr. We maintain the notations introduced above.

THEOREM 10 [Tib02a, Theorem 1.1]. *Let $X_c$ be a general member of the pencil on $X$. Suppose that*

(i) *the axis $A$ is not contained in $V$,*

(ii) *the rectified homotopical depth of $X$ is $\geqslant n$ for some $n \geqslant 2$,*

(iii) *the pencil $\mathscr{P}_A$ has only isolated singular points with respect to $\mathscr{S}$, and*

(iv) *the pair $(X_c, X_c \cap A)$ is $(n-2)$-connected.*

*Then the pair $(X, X_c)$ is $(n-1)$-connected.*

The *rectified homotopical depth* of $X$ is an integer which measures the local connectedness $X$ [HL90, Definition 1.1]. If $X$ is locally a complete intersection variety, it is equal to the complex dimension of $X$ [HL90, Corollary 3.2.2].

**2.3**

Let $S$ be a smooth and irreducible algebraic subset of $\mathbb{P}^n$, and let $A$ be a codimension 2 linear subspace of $\mathbb{P}^n$. We write $\mathscr{P}_A$ for the pencil of hyperplanes containing $A$, and

$$\pi_A : S \setminus A \longrightarrow \mathscr{P}_A$$

for the map sending $p$ to the member of $\mathscr{P}_A$ containing $p$.

LEMMA 11. *If $A_x$ is a general codimension 2 linear subspace passing through a point $x$ in $\mathbb{P}^n$, then*

$$\pi_{A_x} : S \setminus A_x \longrightarrow \mathscr{P}_{A_x}$$

*has only isolated singular points, unless the closure of $S$ in $\mathbb{P}^n$ is a cone with the apex $x$.*

Note that $\pi_{A_x}$ necessarily has nonisolated singularities if, for example, the closure of $S$ in $\mathbb{P}^n$ is the cone over a smooth hypersurface in $\mathbb{P}^{n-1} \subseteq \mathbb{P}^n$ with the apex $x$.

*Proof.* Let $V$ be the closure of $S$ in $\mathbb{P}^n$, and let $A$ be a codimension 2 linear subspace of $\mathbb{P}^n$ containing $x$. Denote the conormal variety of $V$ by $I$, the dual variety of $V$ by $\check{V}$. Consider the projections $\mathrm{pr}_1, \mathrm{pr}_2$ from $I$ to $V, \check{V}$ respectively:

$$\begin{array}{ccccccc}
& & I & & & & \\
& \mathrm{pr}_1 \swarrow & & \searrow \mathrm{pr}_2 & & & \\
V & & & & \check{V} \longrightarrow \check{\mathbb{P}}^n & \longleftarrow & \mathscr{P}_A
\end{array}$$

Choosing a line in $\mathbb{P}^n$ disjoint from $A$ identifies $\pi_A$ with the projection from $A$ to the chosen line. Therefore the singular points of $\pi_A$ are precisely those points at which the projective tangent space of $S$ is contained in a member of the pencil $\mathscr{P}_A$. In other words,

$$\{\text{singular points of } \pi_A\} = \mathrm{pr}_1\big(\mathrm{pr}_2^{-1}(\mathscr{P}_A \cap \check{V})\big) \cap S.$$

Suppose from now on that $V$ is not a cone with the apex $x$. Equivalently, we assume that $\check{V}$ is not contained in $\check{x}$, where $\check{x}$ is a hyperplane in $\check{\mathbb{P}}^n$ corresponding to $x$ [GKZ94, Proposition 4.4].

First consider the case when $\check{V}$ is not a hypersurface in $\check{\mathbb{P}}^n$. In this case, since $\check{V}$ is irreducible and not contained in $\check{x}$,
$$\dim\left(\check{V} \cap \check{x}\right) \leqslant n - 3.$$
Therefore a general line contained in $\check{x}$ is disjoint from $\check{V}$. In other words, $\pi_A$ has no singular points for a general $A$ containing $x$.

Next consider the case when $\check{V}$ is a hypersurface in $\check{\mathbb{P}}^n$. In this case, the biduality theorem shows that $\mathrm{pr}_2$ is generically a projective bundle with zero dimensional fibers [GKZ94, Theorem 1.1]. Since $I$ is irreducible, the previous sentence implies that
$$\dim\left(D := \{y \in \check{V} \mid \mathrm{pr}_2 \text{ has positive dimensional fiber over } y\}\right) \leqslant n - 3.$$
Therefore a general line contained in $\check{x}$ is disjoint from $D$. In other words, $\pi_A$ has only isolated singular points for a general $A$ containing $x$. □

### 2.4

Let $S$ be a smooth and irreducible algebraic subset of $\mathbb{P}^n$, and let $k$ be a positive integer.

LEMMA 12. *If $E_x$ is a general linear subspace of codimension $k$ passing through a point $x$ in $\mathbb{P}^n$, then $E_x$ intersects $S \setminus \{x\}$ transversely in $\mathbb{P}^n$.*

*Proof.* Repeated application of Bertini's theorem shows that $E_x \cap S$ is smooth outside $x$ [Kle98, Theorem 4.1]. In other words, for any $p$ in $E_x \cap S$ different from $x$,
$$\mathrm{codim}\big(\mathrm{T}_p E_x \cap \mathrm{T}_p S \subseteq \mathrm{T}_p S\big) = \mathrm{codim}\big(\mathrm{T}_p E_x \subseteq \mathrm{T}_p \mathbb{P}^n\big) = k.$$
The conclusion follows from the isomorphism
$$\mathrm{T}_p S / \big(\mathrm{T}_p E_x \cap \mathrm{T}_p S\big) \simeq \big(\mathrm{T}_p E_x + \mathrm{T}_p S\big) / \mathrm{T}_p E_x.$$
□

### 2.5

We employ the notation introduced in 2.2. Set $Y = \mathbb{P}^n$, $V = V(h)$, $X = D(h)$, and suppose that

- no component of $V$ is a cone over a smooth variety with the apex $x$, and
- the singular locus of $V$ does not contain a line passing through $x$.

Then we can find a Whitney stratification $\mathscr{W}$ of $Y$ such that

- $\{x\}$ is a stratum of $\mathscr{W}$,
- $V$ is a union of strata of $\mathscr{W}$, and
- the closure of a stratum of $\mathscr{W} \setminus \{\{x\}\}$ is not a cone with the apex $x$.

Let $A$ be a codimension 2 linear subspace of $Y$ containing $x$, and let $\mathbb{Y} \subseteq Y \times \mathbb{P}^1$ be the blow-up of $Y$ along $A$. The projection from $\mathbb{Y}$ onto the $\mathbb{P}^1$ can be identified with the map
$$p : \mathbb{Y} \longrightarrow \mathscr{P}_A.$$
The statement below follows, and can be replaced by, the proof of [Tib02a, Proposition 2.4].

LEMMA 13. *Let $\mathscr{S}$ be the stratification of $\mathbb{Y}$ with strata*

(1) $(S \times \mathbb{P}^1) \cap (\mathbb{Y} \setminus A \times \mathbb{P}^1)$ for $S \in \mathscr{W} \setminus \{\{x\}\}$,

(2) $(S \times \mathbb{P}^1) \cap (A \times \mathbb{P}^1)$ for $S \in \mathscr{W} \setminus \{\{x\}\}$,

(3) $\{x\} \times \mathbb{P}^1 \setminus E$ and $E$,

where $E$ is the set of points at which one of the strata from (1) and (2) fails to be Whitney regular over $\{x\} \times \mathbb{P}^1$. If Lemma 11 and Lemma 12 holds for $A$ and each stratum of $\mathscr{W}$, then

(i) $\mathscr{S}$ is a Whitney stratification, and

(ii) $\mathscr{P}_A$ has only isolated singular points with respect to $\mathscr{S}$.

*Proof.* Let $\mathscr{S}_1$ be the Whitney stratification of $\mathbb{Y} \setminus \{x\} \times \mathbb{P}^1$ with strata

$$\mathbb{Y} \setminus A \times \mathbb{P}^1 \text{ and } A \times \mathbb{P}^1 \setminus \{x\} \times \mathbb{P}^1,$$

and let $\mathscr{S}_2$ be the product Whitney stratification of $Y \times \mathbb{P}^1 \setminus \{x\} \times \mathbb{P}^1$ with strata

$$S \times \mathbb{P}^1 \text{ for } S \in \mathscr{W} \setminus \{\{x\}\}.$$

Lemma 12 shows that any pair of strata from $\mathscr{S}_1$ and $\mathscr{S}_2$ intersect transversely in $Y \times \mathbb{P}^1$. It follows that $\mathscr{S}_1 \cap \mathscr{S}_2$ is a Whitney stratification of $\mathbb{Y} \setminus \{x\} \times \mathbb{P}^1$ [GWPL76, 1.1.3]. Now Whitney's fundamental lemma says that $E$ is finite, and all the strata of $\mathscr{S}_1 \cap \mathscr{S}_2$ are Whitney regular over $E$ [Whi65, Lemma 19.3]. This proves that $\mathscr{S}$ is a Whitney stratification of $\mathbb{Y}$.

For a stratum $\mathcal{S} \in \mathscr{S}$ of the first type, $p|_\mathcal{S}$ has only isolated singular points by Lemma 11. For a stratum $\mathcal{S} \in \mathscr{S}$ of the second type, $p|_\mathcal{S}$ is clearly a submersion and has no singular points, and the same is true for the stratum $\{x\} \times \mathbb{P}^1 \setminus E$. Therefore $\mathscr{P}_A$ has only isolated singular points with respect to $\mathscr{S}$. □

**2.6**

As a final preparation for the proof of Theorem 5, we recall a Zariski theorem on the fundamental group of plane curve complements [Dim01, Section 4.3].

Let $x$ be a point in $\mathbb{P}^2$, and let $C$ be a curve in $\mathbb{P}^2$. We say that a line $L_x$ passing through $x$ is *exceptional* with respect to $C$ if

- $L_x$ is tangent to the curve $C$, or
- $L_x$ passes through a singular point of the curve $C$ different from $x$.

THEOREM 14 [Dim01, Corollary 4.3.6]. *Suppose that no line containing $x$ is a component of the curve $C$. Then for any line $L_x$ passing through $x$ which is not exceptional, there is an epimorphism*

$$\pi_1(L_x \setminus C) \longrightarrow \pi_1(\mathbb{P}^2 \setminus C)$$

*induced by the inclusion.*

Theorem 14 is the base case of the induction for Theorem 5.

**2.7**

*Proof of Theorem 5.* We prove by induction on $n$, the base case being Theorem 14. Suppose that

(a) no component of $V(h)$ is a cone with the apex $x$, and

(b) the singular locus of $V(h)$ does not contain a line passing through $x$.

For the induction step we check that the two conditions on $V(h)$ are also satisfied by $V(h) \cap H_x$. For condition (a), this is the content of Lemma 8. Condition (b) follows from Bertini's theorem that
$$\mathrm{Sing}\big(V(h) \cap H_x\big) \setminus \{x\} = \Big(\mathrm{Sing}\big(V(h)\big) \cap H_x\Big) \setminus \{x\}.$$

Now consider the Whitney stratifications $\mathscr{W}$ and $\mathscr{S}$ of Section 2.5. We also employ other notations introduced in that section. When $n \geqslant 3$, we choose linear subspaces $A \subseteq H$ containing $x$, of codimension 2 and 1 respectively, sufficiently general so that

  (i) $A$ is not contained in $V$,
 (ii) (a) and (b) are satisfied by $V \cap H$,
(iii) Lemma 11 and Lemma 12 holds for $A$ and each stratum of $\mathscr{W}$, and
 (iv) the induction hypothesis applies to the pair $(X \cap H, X \cap A)$.

Then, by Lemma 13, all the assumptions of Theorem 10 are satisfied. Therefore the pair $(X, X \cap H)$ is $(n-1)$-connected. $\square$

## 2.8

For the interested reader we record here a version of what we proved in the generality of Section 2.2. Let $Y$ be a projective variety of dimension $n \geqslant 2$ in $\mathbb{P}^N$, and let $V$ be a closed algebraic subset of $Y$.

THEOREM 15. *Let $E_x \subseteq F_x$ be a general pair of linear subspaces containing a point $x$ in $\mathbb{P}^N$, of codimensions $n-1$ and $n-2$ respectively. Suppose that*

  (i) *the quasi-projective variety $X := Y \setminus V$ is locally a complete intersection,*
 (ii) *no component of $V$ (and of $Y$) is a cone with the apex $x$,*
(iii) *the singular locus of $V$ (and of $Y$) does not contain a line passing through $x$, and*
 (iv) *there is an epimorphism induced by the inclusion*
$$\pi_1\big(X \cap E_x\big) \longrightarrow \pi_1\big(X \cap F_x\big).$$

*Then, for a sufficiently general hyperplane $H_x$ passing through $x$,*
$$\pi_i\big(X, X \cap H\big) = 0 \quad \text{for} \quad i < n.$$

## 3. Proof of Theorem 2

## 3.1

We deduce Theorem 2 from Theorem 5 when $n \geqslant 3$. A separate argument will be given for plane curves.

*Proof of Theorem 2 when $n \geqslant 3$.* We know from [DP03] that
$$\chi\big(D(h)\big) = (-1)^n \deg\big(\mathrm{grad}(h)\big) + \sum_{i=0}^{n-1} (-1)^i (d-1)^i.$$

If $V(h)$ is not a cone with the apex $x$, then there is a hyperplane $H_x$ containing $x$ such that

  (i) Theorem 5 applies to $H_x$,

(ii) $V(h) \cap H_x$ is smooth outside $x$, and

(iii) the Milnor number of $V(h) \cap H_x$ at $x$ is the sectional Milnor number $\mu^{(n-1)}$ of $V(h)$ at $x$.

It follows from (ii) and (iii) that

$$\chi\big(D(h) \cap H_x\big) = (-1)^{n-1}\Big((d-1)^{n-1} - \mu^{(n-1)}\Big) + \sum_{i=0}^{n-2}(-1)^i(d-1)^i.$$

Therefore

$$\text{rank } H_n\big(D(h), D(h) \cap H_x\big) = (-1)^n\Big(\chi\big(D(h)\big) - \chi\big(D(h) \cap H_x\big)\Big) = \deg\big(\text{grad}(h)\big) - \mu^{(n-1)} \geqslant 0.$$

$\square$

**3.2**

We prove Theorem 2 for plane curves. The main ingredient in this case is Milnor's formula for the double point number

$$2\delta_x = \mu_x + r_x - 1,$$

where $\mu_x$ is the Milnor number at $x$, and $r_x$ is the number of branches at $x$ [Mil68, Theorem 10.5].

LEMMA 16.

(i) Suppose $V(h)$ is a reduced and irreducible plane curve of degree $d$ containing $x$. Then

$$\deg\big(\text{grad}(h)\big) \geqslant (d-1) + (r_x - 1),$$

where $r_x$ is the number of branches of $V(h)$ at $x$.

(ii) Suppose $V(h_1)$ and $V(h_2)$ are plane curves with no common components. Then

$$\deg\big(\text{grad}(h_1 h_2)\big) = \deg\big(\text{grad}(h_1)\big) + \deg\big(\text{grad}(h_2)\big) + \#\big(V(h_1) \cap V(h_2)\big) - 1.$$

*Proof.* The first inequality is obtained from [Mil68, Theorem 10.5] and the fact that $g \geqslant 0$, where $g$ is the genus of the normalization of $V(h)$.

The second assertion is equivalent to the inclusion-exclusion formula for the topological Euler characteristic.

We refer to [Dol00, Section 3] and [FM12, Theorem 3.1] for details. $\square$

*Proof of Theorem 2 when $n = 2$.* Let $V$ be a reduced plane curve of degree $d$ which is not a cone with the apex $x$. Lemma 16 proves the assertion when $V$ is irreducible. We divide the remaining problem into two cases:

(1) $V$ has at least two components which are not cones with the apex with $x$.

(2) $V$ has exactly one component which is not a cone with the apex with $x$.

In case (1), we induct on the number of irreducible components. Write $V = V_1 \cup V_2$, where $V_1$ is not a cone with the apex $x$, $V_2$ is irreducible and not a cone with the apex $x$, and $V_1 \cap V_2$ is finite. Let $h_1, h_2$ be reduced equations of degree $d_1, d_2$ defining $V_1, V_2$ respectively. Then, by Lemma 16 and the induction hypothesis,

$$\deg\big(\text{grad}(h_1 h_2)\big) = \deg\big(\text{grad}(h_1)\big) + \deg\big(\text{grad}(h_2)\big) + \Big(\#(V_1 \cap V_2) - 1\Big)$$

$$\geqslant \Big(m_x(V_1) - 1\Big) + \Big(d_2 - 1\Big) + \Big(\#(V_1 \cap V_2) - 1\Big).$$

In other words,
$$\deg\bigl(\mathrm{grad}(h_1 h_2)\bigr) \geqslant \bigl(m_x(V) - 1\bigr) + \bigl(d_2 - m_x(V_2) - 1\bigr) + \bigl(\#(V_1 \cap V_2) - 1\bigr).$$
The second term in the last expression is nonnegative because $V_2$ is not a cone with the apex $x$. The third term is also nonnegative, and this gives the desired inequality.

In case (2), we write $V = V_1 \cup V_2$, where $V_1$ is a cone with the apex $x$, and $V_2$ is irreducible. Then, by Lemma 16,
$$\deg\bigl(\mathrm{grad}(h_1 h_2)\bigr) = \deg\bigl(\mathrm{grad}(h_2)\bigr) + \bigl(\#(V_1 \cap V_2) - 1\bigr)$$
$$\geqslant \bigl(d_2 - 1\bigr) + \bigl(r_x(V_2) - 1\bigr) + \bigl(\#(V_1 \cap V_2) - 1\bigr).$$

In other words,
$$\deg\bigl(\mathrm{grad}(h_1 h_2)\bigr) \geqslant \bigl(m_x(V) - 1\bigr) + \bigl(d_2 - m_x(V_2) - 1\bigr) + \bigl(r_x(V_2) - m_x(V_1) + \#(V_1 \cap V_2) - 1\bigr).$$
The second term in the last expression is nonnegative because $V_2$ is not a cone with the apex $x$. We claim that the third term is also nonnegative.

Let $t_x(V_2)$ be the number of lines in the tangent cone of $V_2$ at $x$. It follows from Hensel's lemma that
$$r_x(V_2) \geqslant t_x(V_2).$$
See for example [Cas00, Corollary 2.2.6]. A local computation shows that a line containing $x$ intersects $V_2$ in at least one point other than $x$, unless the line is contained in the tangent cone of $V_2$ at $x$. Since there are at least $m_x(V_1) - t_x(V_1)$ lines in $V_1$ not contained in the tangent cone of $V_2$ at $x$, we have
$$t_x(V_2) - m_x(V_1) + \#\bigl(V(h_1) \cap V(h_2)\bigr) - 1 \geqslant 0.$$
This completes the proof. □

### 3.3

We show that, for each $d \geqslant 3$ and $n \geqslant 2$, there is a degree $d$ hypersurface in $\mathbb{P}^n$ with one singular point, at which the equality holds in Theorem 2.

PROPOSITION 17. *Let $V(h)$ be the degree $d$ hypersurface in $\mathbb{P}^n$ defined by the equation*
$$h = z_0 z_1^{d-1} + z_1 z_2^{d-1} + (z_3^d + \cdots + z_n^d), \qquad d \geqslant 3.$$
*Then the unique singular point $x$ of $V(h)$ satisfies*
$$\mu^{(n)} = (d-1)^n - (d-1)^{n-1} + (d-1)^{n-2} \quad \text{and} \quad \mu^{(n-1)} = (d-1)^{n-1} - (d-1)^{n-2}.$$

*Proof of Proposition 17 when $n \geqslant 3$.* Locally at $x$, the hypersurface is defined by
$$f = x_1^{d-1} + x_1 x_2^{d-1} + x_3^d + \cdots + x_n^d.$$
Note that $f$ is weighted homogeneous with weights $w_i$, where
$$\frac{d-1}{w_1} = \frac{1}{w_1} + \frac{d-1}{w_2} = \frac{d}{w_3} = \cdots = \frac{d}{w_n} = 1.$$
It follows from [MO70, Theorem 1] that
$$\mu^{(n)} = \prod_{i=1}^{n}(w_i - 1) = (d-1)^n - (d-1)^{n-1} + (d-1)^{n-2}.$$

Now consider the hyperplane $H$ passing through $x$ defined by

$$x_n = c_1 x_1 + \cdots + c_{n-1} x_{n-1}, \qquad c = (c_1, \ldots, c_{n-1}) \in (\mathbb{C}^*)^{n-1}.$$

Locally at $x$, $V(h) \cap H$ is isomorphic to the hypersurface defined by

$$g = x_1^{d-1} + x_1 x_2^{d-1} + x_3^d + \cdots + x_{n-1}^d + (c_1 x_1 + \cdots + c_{n-1} x_{n-1})^d.$$

The principal part of $g$ with respect to the Newton diagram is

$$g_0 = x_1^{d-1} + x_3^d + \cdots + x_{n-1}^d + (c_2 x_2 + \cdots + c_{n-1} x_{n-1})^d.$$

Since $g_0$ defines an isolated singular point at the origin, $g$ is semiquasihomogeneous. This shows that the singular points defined by $g$ and $g_0$ have the same Milnor number

$$\mu^{(n-1)} = (d-1)^{n-1} - (d-1)^{n-2}.$$

See [AGV85, Chapter 12]. □

Proposition 17 remains valid for $n = 2$.

## 4. Projective hypersurfaces with small polar degree

### 4.1

A reduced homogeneous polynomial $h$ is homaloidal if and only if $n$ general first polar hypersurfaces of $V(h)$ intersect at exactly one point outside the singular locus of $V(h)$. However, Conjecture 3 cannot be proved by a Cayley-Bacharach type theorem alone.

*Example* 18. For each $d \geqslant 2$ and $n \geqslant 2$, there is a birational transformation

$$\varphi : \mathbb{P}^n \dashrightarrow \mathbb{P}^n, \qquad z \longmapsto (h_0 : \cdots : h_n),$$

where $h_0, \ldots, h_1$ are homogeneous polynomials of the same degree $d-1$, with a zero-dimensional base locus. For example, one may take

$$h_n = z_0^{d-1}, \qquad h_{n-1} = 2 z_0^{d-2} z_1, \qquad h_{n-i} = 2 z_0^{d-2} z_i + z_{i-1}^{d-1}, \qquad i = 2, \ldots, n.$$

When $n = 2$ and $d = 3$, this is the gradient map of the homaloidal polynomial

$$h = z_0(z_1^2 + z_0 z_2).$$

Conjecture 3 asserts that no birational transformation with a zero-dimensional base locus is a gradient map when $d \geqslant 3$ and $n \geqslant 3$.

### 4.2

A plane curve singularity of multiplicity 2 is locally defined by an equation of the form

$$f = x_1^2 + x_2^{k+1},$$

where $k$ is the Milnor number at the singular point [Har77, Exercise 1.5.14]. Similar statements remain valid in all dimensions.

LEMMA 19. *Let $(V, \mathbf{0})$ be the germ of an isolated hypersurface singularity at the origin of $\mathbb{C}^n$.*
  (i) *If $\mu^{(n-1)}$ of the germ is equal to 1, then the singularity is of type $A_k$ for some $k \geqslant 1$.*
  (ii) *If $\mu^{(n-1)}$ of the germ is equal to 2, then the singularity is of type $E_{6r}$, $E_{6r+1}$, $E_{6r+2}$ or $J_{r,i}$, for some $r \geqslant 1$ and $i \geqslant 0$.*

For the normal form of the above singularities, we refer to [AGV85, Chapter 15]. We set $J_{1,i} = D_{4+i}$ in the second case so that the singularity is simple if and only if $r = 1$.

*Proof.* We prove the first statement. The second statement can be justified in the same way.

Choose any hyperplane $H$ passing through the origin such that $(V \cap H, \mathbf{0})$ has the Milnor number 1. Let $f$ be an equation defining $V$, and let $y_n$ be an equation defining $H$. Morse lemma shows that there is a local coordinate system $y_1, \ldots, y_{n-1}$ of $(H, \mathbf{0})$ such that

$$f|_H = y_1^2 + \cdots + y_{n-2}^2 + y_{n-1}^2.$$

Therefore, we may write

$$f = y_1^2 + \cdots + y_{n-1}^2 + y_n(c_1 y_1 + \cdots + c_n y_n) + \big(\text{terms of degree} \geqslant 3 \text{ in } y_1, \ldots, y_n\big)$$

for some $c_1, \ldots, c_n \in \mathbb{C}$. This shows that the Hessian of $f$ at $\mathbf{0}$ relative to $y_1, \ldots, y_n$ is

$$H(f) = \begin{pmatrix} 1 & 0 & \cdots & 0 & c_1 \\ 0 & 1 & \cdots & 0 & c_2 \\ \vdots & \vdots & & \vdots & \vdots \\ 0 & 0 & \cdots & 1 & c_{n-1} \\ c_1 & c_2 & \cdots & c_{n-1} & c_n \end{pmatrix}.$$

In particular, the corank of the Hessian is at most 1. Now the classification of corank $\leqslant 1$ isolated singularities says that there is a local coordinate system $x_1, \ldots, x_n$ of $(\mathbb{C}^n, \mathbf{0})$ with

$$f = x_1^2 + \cdots + x_{n-1}^2 + x_n^{k+1},$$

where $k$ is the Milnor number of $V$ at the origin. See [AGV85, Chapter 11]. □

### 4.3

Dimca proves Conjecture 3 in [Dim01, Theorem 9] when all the singular points of $V(h)$ are weighted homogeneous. The proof is based on the work of du Plessis and Wall [dPW01]. We use Theorem 2 to reduce the problem to this case.

*Proof of Theorem 4.* Let $V(h)$ be a homaloidal hypersurface of degree $d$ in $\mathbb{P}^n$ with only isolated singular points. Theorem 2 implies that all the singular points of $V(h)$ has the sectional Milnor number $\mu^{(n-1)} = 1$. It follows from Lemma 19 that all the singular points of $V(h)$ are locally defined by an equation of the form

$$f = x_1^2 + \cdots + x_{n-1}^2 + x_n^{k+1},$$

where $k$ is the Milnor number of $V(h)$ at the singular point. In particular, all the singular points of $V(h)$ are weighted homogeneous. Therefore, by [Dim01, Theorem 9], either $n \leqslant 2$ or $d \leqslant 2$. If $d \leqslant 2$, the hypersurface should be smooth quadric, and any smooth quadric is defined by

$$h = z_0^2 + \cdots + z_n^2,$$

after a linear change of coordinates. When $n \leqslant 2$, the assertion is [Dol00, Theorem 4]. □

### 4.4

What can be said about projective hypersurfaces which has polar degree 2 and only isolated singular points? We propose the following conjecture.

CONJECTURE 20. *A projective hypersurface with only isolated singular points has polar degree 2 if and only if it is one of the following, after a linear change of homogeneous coordinates:*

- ($n = 3, d = 3$) *a normal cubic surface containing a single line*
$$h = z_0 z_1^2 + z_1 z_2^2 + z_1 z_3^2 + z_2^3 = 0, \qquad (E_6).$$

- ($n = 3, d = 3$) *a normal cubic surface containing two lines*
$$h = z_0 z_1 z_2 + z_0 z_3^2 + z_1^3 = 0, \qquad (A_5, A_1).$$

- ($n = 3, d = 3$) *a normal cubic surface containing three lines and three binodes*
$$h = z_0 z_1 z_2 + z_3^3 = 0, \qquad (A_2, A_2, A_2).$$

- ($n = 2, d = 5$) *two smooth conics meeting at a single point and the common tangent*
$$h = z_0(z_1^2 + z_0 z_2)(z_1^2 + z_0 z_2 + z_0^2) = 0, \qquad (J_{2,4}).$$

- ($n = 2, d = 4$) *two smooth conics meeting at a single point*
$$h = (z_1^2 + z_0 z_2)(z_1^2 + z_0 z_2 + z_0^2) = 0, \qquad (A_7).$$

- ($n = 2, d = 4$) *a smooth conic, a tangent, and a line passing through the tangency point*
$$h = z_0(z_0 + z_1)(z_1^2 + z_0 z_2) = 0, \qquad (D_6, A_1).$$

- ($n = 2, d = 4$) *a smooth conic and two tangent lines*
$$h = z_0 z_2 (z_1^2 + z_0 z_2) = 0, \qquad (A_1, A_3, A_3).$$

- ($n = 2, d = 4$) *three concurrent lines and a line not meeting the center point*
$$h = z_0 z_1 z_2 (z_0 + z_1) = 0, \qquad (D_4, A_1, A_1, A_1).$$

- ($n = 2, d = 4$) *a cuspidal cubic and its tangent at the cusp*
$$h = z_0(z_1^3 + z_0^2 z_2) = 0, \qquad (E_7).$$

- ($n = 2, d = 4$) *a cuspidal cubic and its tangent at the smooth flex point*
$$h = z_2(z_1^3 + z_0^2 z_2) = 0, \qquad (A_2, A_5).$$

- ($n = 2, d = 3$) *a cuspidal cubic*
$$h = z_1^3 + z_0^2 z_2 = 0, \qquad (A_2).$$

- ($n = 2, d = 3$) *a smooth conic and a secant line*
$$h = z_1(z_1^2 + z_0 z_2) = 0, \qquad (A_1, A_1).$$

Less precisely but more generally, we conjecture that for any positive integer $k$, there is no projective hypersurface of polar degree $k$ which has only isolated singular points, for sufficiently large $n$ and $d$.

PROPOSITION 21. *Conjecture 20 is valid for plane curves, cubic surfaces, and quartic surfaces.*

Radu Laza informed us that Conjecture 20 is valid also for cubic threefolds. We note that there is a cubic threefold with only isolated singular points which has polar degree 3. Explicitly, there is
$$h = z_0 z_1 z_4 + z_0^3 + z_1^3 + z_0 z_2^2 + z_1 z_3^2 = 0, \qquad (T_{2,6,6}).$$

*Proof.* The proof is a combination of [BW79, Deg90, Dim01, dPW01, FM12, Wal99].

The assertion for cubic surfaces is classical, and can be deduced from the classification [BW79]. The list of plane curves is obtained in [FM12]. Any reduced plane curve with polar degree 2 should be projectively equivalent to one of the list.

For the remaining case of quartic surfaces, we start with two general remarks. We assume throughout that projective hypersurfaces have only isolated singular points.

(i) If $V(h)$ is a projective hypersurface with only weighted homogeneous singular points, then
$$\deg\bigl(\mathrm{grad}(h)\bigr) \geqslant \min\Bigl\{(d-1)^{n-2}, 2(n+1)\Bigr\},$$
unless $V(h)$ is a cone. This follows from the proof of [Dim01, Theorem 9]. We note from the above example of cubic threefold with $T_{2,6,6}$ singularity that the assumption of weighted homogeneity is necessary for the inequality.

(ii) Each nonsimple singular point of a projective hypersurface of polar degree 2 should belong to one of the types $E_{6r}$, $E_{6r+1}$, $E_{6r+2}$ or $J_{r,i}$ for some $r \geqslant 2$ and $i \geqslant 0$. This follows from Theorem 2 and Lemma 19.

Let $V$ be a quartic surface of polar degree 2. The first remark above shows that it is enough to consider the case when $V$ has a singular point $p$ which is not simple. The second remark shows that the singularity of $V$ at $p$ is of type **J** or **E**. Let $C$ be the discriminant of the projection $V \setminus \{p\} \to \mathbb{P}^2$ from $p$. Under our assumptions, $C$ is a reduced sextic curve.

Suppose $V$ is not stable in the sense of Geometric Invariant Theory. In this case, there is a bijection between the singular points of $V$ and of $C$, which preserves the type [Wal99, Section 11]. In particular,
$$\sum_{x \in V} \mu^{(3)}(x, V) = \sum_{x \in C} \mu^{(2)}(x, C).$$
Since $V$ has polar degree 2, the sum of the Milnor numbers of $C$ should be $3^3 - 2 = 5^2$. From Theorem 1 it follows that $C$ is a cone. This contradicts the fact that $C$ has a singular point of type **J** or **E**. In Degtyarëv's classification of quartic surfaces having a nonsimple singular points, this case is called *exceptional* [Deg90, Theorem 1.9].

Suppose $V$ is stable in the sense of Geometric Invariant Theory. In this case, there is a point $q$ such that the blowup of the singularity of $V \subseteq \mathbb{P}^3$ at $p$ has the same type as the singularity of $C \subseteq \mathbb{P}^2$ at $q$. Moreover, there is a bijection between the singular points of $V \setminus \{p\}$ and of $C \setminus \{q\}$, which preserves the type [Wal99, Section 12]. A case by case analysis of the blowup of suspensions of singularities of type **J** and **E** shows that
$$\sum_{x \in V} \mu^{(3)}(x, V) = 1 + \sum_{x \in C} \mu^{(2)}(x, C).$$
Since $V$ has polar degree 2, the sum of the Milnor numbers of $C$ should be $3^3 - 3 = 5^2 - 1$. This means that the sextic plane curve $C$ is homaloidal, which is impossible. In Degtyarëv's classification of quartic surfaces having a nonsimple singular points, this case is called *nonexceptional* [Deg90, Theorem 1.7]. □

June Huh   junehuh@umich.edu

Department of Mathematics, University of Michigan, Ann Arbor, MI 48109, USA